\documentclass[11pt,reqno]{amsart}
\usepackage{amsmath,bm} 

\usepackage{amssymb,latexsym}
\usepackage{graphicx}
\usepackage{fancyhdr,amssymb}
\usepackage{amstext}
\usepackage{url}		

\theoremstyle{plain}
\numberwithin{equation}{section}
\newtheorem{thm}{Theorem}[section]
\newtheorem{theorem}[thm]{Theorem}
\newtheorem{lemma}[thm]{Lemma}

\newtheorem{corollary}[thm]{Corollary}
\newtheorem{remark}[thm]{Remark}
\newtheorem{conjecture}[thm]{Conjecture}

\def\dire{dir{\hspace{-0.11em}\raisebox{-0.65ex}{\mbox{\scriptsize $\vec{e}$}}}\,}

\begin{document}
\fancyhead{}
\renewcommand{\headrulewidth}{0pt}
\fancyfoot{}
\fancyfoot[LE,RO]{\medskip \thepage}
\fancyfoot[LO]{\medskip }
\fancyfoot[RE]{\medskip }

\setcounter{page}{1}

\title[SOME ENUMERATIONS OF NON-TRIVIAL COMPOSITIONS OF THE DIFFERENTIAL OPERATIONS AND THE DIRECTIONAL DERIVATIVE]{
 \mbox{SOME ENUMERATIONS OF NON-TRIVIAL COMPOSITIONS OF}                        \\[0.05 ex]
 \mbox{THE DIFFERENTIAL OPERATIONS AND}                                         \\[0.05 ex]
 \mbox{THE~DIRECTIONAL~DERIVATIVE}}
\author{Ivana Jovovi\' c}
\address{University of Belgrade, Faculty of Electrical Engineering,             \newline
         \hspace*{8.0 mm}Bulevar Kralja Aleksandra 73, 11000 Belgrade, Serbia}
\email{ivana.jovovic@etf.bg.ac.rs}
\author{Branko Male\v sevi\' c}
\address{University of Belgrade, Faculty of Electrical Engineering,             \newline
         \hspace*{8.0 mm}Bulevar Kralja Aleksandra 73, 11000 Belgrade, Serbia}
\email{branko.malesevic@etf.bg.ac.rs}

\begin{abstract}
This paper deals with some enumerations of the higher order non-trivial compositions of
the differential operations and the directional derivative in the space $\mathbb{R}^n$
($n \!\ge\! 3$). One~new~enumeration of the higher order non-trivial compositions is obtained.
\end{abstract}

\noindent
{\small \tt Notes on Number Theory}

\vspace*{-0.25 mm}

\noindent
{\small \tt \& Discrete Mathematics \\ } 

\bigskip
\maketitle

\vspace*{-5.0 mm}

\section{The enumeration of the higher order non-trivial compositions of the differential operations
         and the directional derivative in the space $\mathbb{R}^3$}

\smallskip
Consider the sets of the smooth functions
$$
\mbox{\rm A}_{0}
\!=\!
\bigl\{
f\!:\! \mathbb{R}^3 \!\rightarrow\! \mathbb{R}
\, | f\!\in\!C^{\infty}(\mathbb{R}^3)
\bigr\}
\;\;
\mbox{and}
\;\;
\mbox{\rm A}_{1}
\!=\!
\bigl\{
\vec{f}\!:\! \mathbb{R}^3 \!\rightarrow\! \mathbb{R}^3
\, | f_{1},f_{2},f_{3}\!\in\!C^{\infty}(\mathbb{R}^3)
\bigr\}
$$
in the three-dimensional Euclidean space $\mathbb{R}^{3}$.
Let $\vec{e} = (e_1,e_2,e_3) \in \mathbb{R}^{3}$ be a unit vector.
The gradient, curl, divergence, and the Gateaux directional derivative in a direction $\vec{e}$
are defined in the terms of the partial derivative operators as follows:
$$
\begin{array}{lclcll}
\mbox{\rm grad} \, f
\!\!\!&\!=\!&\!\!\!
\nabla_1 \, f
\!\!\!&=&\!\!\!
\frac{\partial f}{\partial x_1} \, \vec{i} \!+\!
\frac{\partial f}{\partial x_2} \, \vec{j} \!+\!
\frac{\partial f}{\partial x_3} \, \vec{k},   \;
\nabla_1 \!:\! \mbox{\rm A}_{0} \!\longrightarrow\! \mbox{\rm A}_{1} \, ,                             \\[1.25 ex]
\!\!\mbox{ \rm curl} \, \vec{f}
\!\!&\!=\!&\!\!\!
\nabla_2 \, \vec{f}
\!\!\!&=&\!\!\!
\left(\! \frac{\partial f_3}{\partial x_2} \!-\!
\frac{\partial f_2}{\partial x_3} \!\right) \vec{i} \!+\!
\left(\! \frac{\partial f_1}{\partial x_3} \!-\!
\frac{\partial f_3}{\partial x_1} \!\right) \vec{j}  \!+\!
\left(\! \frac{\partial f_2}{\partial x_1} \!-\!
\frac{\partial f_1}{\partial x_2} \!\right) \vec{k},    \;
\nabla_2 \!:\! \mbox{\rm A}_{1} \!\longrightarrow\! \mbox{\rm A}_{1} \, ,                             \\[1.25 ex]
\mbox{\rm div} \, \vec{f}
\!\!\!&\!=\!&\!\!\!
\nabla_3 \, \vec{f}
\!\!\!&=&\!\!\!
\frac{\partial f_1}{\partial x_1} \!+\!
\frac{\partial f_2}{\partial x_2} \!+\!
\frac{\partial f_3}{\partial x_3},   \;
\nabla_3 \!:\! \mbox{\rm A}_{1} \!\longrightarrow\! \mbox{\rm A}_{0} \, ,                             \\[1.25 ex]
\mbox{\dire} \, f
\!\!\!&\!=\!&\!\!\!
\nabla_0 \, f
\!\!\!&=&\!\!\!
 \nabla_1 \, f \cdot \vec{e} =
\frac{\partial f}{\partial x_1} \, e_1 \!+\!
\frac{\partial f}{\partial x_2} \, e_2 \!+\!
\frac{\partial f}{\partial x_3} \, e_3,   \;
\nabla_0 \!:\! \mbox{\rm A}_{0} \!\longrightarrow\! \mbox{\rm A}_{0} \, .
\end{array}
$$

Let $\mathcal{A}_{3} \!=\! \{ \nabla_1, \nabla_2, \nabla_3 \}$ and
$\mathcal{B}_{3} \!=\! \{ \nabla_0, \nabla_1, \nabla_2, \nabla_3 \}$.
Male\v sevi\' c~\cite{HiOrd98} has proved that the number of the
$k^{\mbox{\scriptsize \rm th}}$ order compositions over the set
$\mathcal{A}_{3}$ is $\mbox{\textbf{\texttt{f}}}(k) = F_{k+3}$,
$F_{k}$ is the $k^{\mbox{\scriptsize \rm th}}$ Fibonacci number.
A composition of differential operations that is not equal to $0$ or $\vec{0}$ is called non-trivial.
Male\v sevi\' c~\cite{HiOrd96} has showed that
the number of the $k^{\mbox{\scriptsize \rm th}}$ order non-trivial compositions
over the set $\mathcal{A}_{3}$ is $\mbox{\textbf{\texttt{g}}}(k) = 3$.
Schreiber \cite[Section 5.2]{Schreiber77} has listed
the higher order non-trivial compositions over the set $\mathcal{A}_{3}$:
$$
\begin{array}{rcl}
\mbox{(grad)\,div\,\ldots\,grad\,div\,grad}\,     f \!&\!=\!&\! (\nabla_1\circ) \nabla_3 \circ \cdots \circ \nabla_1 \circ \nabla_3 \circ \nabla_1 f,      \\[1.0 ex]
\mbox{curl\,curl\,\ldots\,curl\,curl\,curl}\,\vec f \!&\!=\!&\! \nabla_2\circ \nabla_2 \circ \cdots \circ \nabla_2\circ\nabla_2 \circ \nabla_2 \vec f,     \\[1.0 ex]
\mbox{(div)\,grad\,\ldots\,div\,grad\,div}\,\vec  f \!&\!=\!&\! (\nabla_3\circ) \nabla_1 \circ \cdots \circ \nabla_3\circ\nabla_1 \circ \nabla_3 \vec f.
\end{array}
$$
The terms in brackets are included if the number of the differential operations is odd and are left out otherwise.
Male\v sevi\' c, and Jovovi\' c \cite{HiOrd07} have proved that the number of the $k^{\mbox{\scriptsize \rm th}}$ order compositions
over the set $\mathcal{B}_{3}$ is $\mbox{\textbf{\texttt{f}}}^{\mbox{\tiny \rm G}\!}(k) = 2^{k+1}$.

\smallskip

According to the above results, it is natural to try to calculate the number of the non-trivial compositions over the set $\mathcal{B}_{3}$.
A straightforward verification shows that all compositions of the second order over $\mathcal{B}_{3}$ are
$$
\begin{array}{l}
\mbox{\dire\,\dire} \, f = \nabla_0 \circ \nabla_0\, f
=
\nabla_1 {\big (} \, \nabla_1 f \cdot \vec{e} \, {\big )} \cdot \vec{e},                          \\[1.0 ex]
\mbox{grad\,\dire} \, f = \nabla_1 \circ \nabla_0\, f
=
\nabla_1 {\big (} \, \nabla_1 f \cdot \vec{e} \, {\big )},                                        \\[1.0 ex]
\mbox{\dire\,div} \, \vec{f} = \nabla_0 \circ \nabla_3\, \vec{f}
=
{\big (} \nabla_1 \circ \nabla_3 \, \vec{f} \, {\big )} \cdot \vec{e},                            \\[1.0 ex]
\mbox{grad\,div} \, \vec{f} = \nabla_1 \circ \nabla_3\, \vec{f},                                  \\[1.0 ex]
\mbox{curl\,curl} \, \vec{f} = \nabla_2 \circ \nabla_2\, \vec{f},                                 \\[1.0 ex]
\mbox{div\,grad} \, f = \nabla_3 \circ \nabla_1 \, f = \Delta f,                                  \\[1.0 ex]
\mbox{curl\,grad} \, f = \nabla_2 \circ \nabla_1\, f = \vec{0},                                   \\[1.0 ex]
\mbox{div\,curl} \, \vec{f} = \nabla_3 \circ \nabla_2\, \vec{f} = 0,
\end{array}
$$
and that only the last two are trivial.
This fact leads us to use the following method for determining the number of
the non-trivial compositions over the set $\mathcal{B}_{3}$.
We define a binary relation $\sigma$ on the set $\mathcal{B}_{3}$ as follows:
\begin{center}
$\nabla_{i} \,\sigma\, \nabla_{j}$ iff the composition $\nabla_{j} \circ \nabla_{i}$ is non-trivial.
\end{center}

The relation $\sigma$ induces the Cayley table

\vspace*{-4.0 mm}

$$
\begin{array}{c|cccc}
\sigma & \nabla_{0} \!&\! \nabla_{1} \!&\! \nabla_{2} \!&\! \nabla_{3} \\[0.11 ex] \hline \\[-2.0 ex]
\nabla_{0}          \!&\! 1          \!&\! 1          \!&\! 0          \!&\! 0            \\[0.75 ex]
\nabla_{1}          \!&\! 0          \!&\! 0          \!&\! 0          \!&\! 1            \\[0.75 ex]
\nabla_{2}          \!&\! 0          \!&\! 0          \!&\! 1          \!&\! 0            \\[0.75 ex]
\nabla_{3}          \!&\! 1          \!&\! 1          \!&\! 0          \!&\! 0
\end{array}
$$

For a convenience, we extend the set $\mathcal{B}_{3}$ with the nowhere-defined function $\nabla_{\!\!-1}$, whose
domain and range are empty sets, and establish $\nabla_{\!\!-1} \,\sigma\, \nabla_{i}$ ($0 \leq i \leq 3$).
Thus, the graph $\textmd{$\Gamma$}$ of the relation $\sigma$ is rooted tree with a root $\nabla_{\!\!-1}$.

\setlength{\unitlength}{0.19 cc}                                                             
\begin{picture}(140,48)(0,0)                                                                 
\thicklines                                                                                  
\put(75,40){\circle*{0.8}}                                                                   
\put(76,41){\scriptsize$\nabla_{\!\!-1}$}                                                    
\put(152,40){\scriptsize $\mbox{\small \textbf{\texttt{g}}}^{\mbox{\tiny \rm G}\!}(0)=1$}    
\put(25,30){\line(5,1){50}}                                                                  
\put(25,30){\circle*{0.8}}                                                                   
\put(23,32){\scriptsize$\nabla_{0}$}                                                         
\put(55,30){\line(2,1){20}}                                                                  
\put(55,30){\circle*{0.8}}                                                                   
\put(52,32){\scriptsize$\nabla_{1}$}                                                         
\put(95,30){\line(-2,1){20}}                                                                 
\put(95,30){\circle*{0.8}}                                                                   
\put(95,32){\scriptsize$\nabla_{2}$}                                                         
\put(125,30){\line(-5,1){50}}                                                                
\put(125,30){\circle*{0.8}}                                                                  
\put(125,32){\scriptsize$\nabla_{3}$}                                                        
\put(152,30){\scriptsize $\mbox{\small \textbf{\texttt{g}}}^{\mbox{\tiny \rm G}\!}(1)=4$}    
\put(15,20){\line(1,1){10}}                                                                  
\put(15,20){\circle*{0.8}}                                                                   
\put(11,22){\scriptsize$\nabla_{0}$}                                                         
\put(35,20){\line(-1,1){10}}                                                                 
\put(35,20){\circle*{0.8}}                                                                   
\put(35,22){\scriptsize$\nabla_{1}$}                                                         
\put(55,20){\line(0,1){10}}                                                                  
\put(55,20){\circle*{0.8}}                                                                   
\put(56,22){\scriptsize$\nabla_{3}$}                                                         
\put(95,20){\line(0,1){10}}                                                                  
\put(95,20){\circle*{0.8}}                                                                   
\put(96,22){\scriptsize$\nabla_{2}$}                                                         
\put(115,20){\line(1,1){10}}                                                                 
\put(115,20){\circle*{0.8}}                                                                  
\put(111,22){\scriptsize$\nabla_{0}$}                                                        
\put(135,20){\line(-1,1){10}}                                                                
\put(135,20){\circle*{0.8}}                                                                  
\put(135,22){\scriptsize$\nabla_{1}$}                                                        
\put(152,20){\scriptsize $\mbox{\small \textbf{\texttt{g}}}^{\mbox{\tiny \rm G}\!}(2)=6$}    
\put(12.5,15){\line(1,2){2.5}}                                                               
\put(17.5,15){\line(-1,2){2.5}}                                                              
\put(35,15){\line(0,1){5}}                                                                   
\put(52.5,15){\line(1,2){2.5}}                                                               
\put(57.5,15){\line(-1,2){2.5}}                                                              
\put(95,15){\line(0,1){5}}                                                                   
\put(112.5,15){\line(1,2){2.5}}                                                              
\put(117.5,15){\line(-1,2){2.5}}                                                             
\put(135,15){\line(0,1){5}}                                                                  
\put(152,12){\scriptsize $\mbox{\small \textbf{\texttt{g}}}^{\mbox{\tiny \rm G}\!}(3)=9$}    
\thinlines                                                                                   
\put(65,6.5){\small Fig. 1. Tree $\textmd{$\Gamma$}$}                                        
\end{picture} 

\vspace*{-1.5 mm}

\noindent
Here we would like to point out, that the child of $\nabla_{i}$ is $\nabla_{j}$
if composition $\nabla_{j} \circ \nabla_{i}$ is non-trivial. For any
non-trivial composition $\nabla_{i_k} \circ \cdots \circ \nabla_{i_1}$ there is
a unique path in the tree $\textmd{$\Gamma$}$, such that the level of vertex $\nabla_{i_j}$ is $j$
($1\leq j \leq k$).

Let $\mbox{\textbf{\texttt{g}}}^{\mbox{\tiny \rm G}\!}(k)$ be the number of the $k^{\mbox{\scriptsize \rm th}}$ order
non-trivial compositions over the set $\mathcal{B}_{3}$.
Let $\mbox{\textbf{\texttt{g}}}_{i}^{\mbox{\tiny \rm G}\!}(k)$ be the number of the $k^{\mbox{\scriptsize \rm th}}$ order
non-trivial compositions starting with $\nabla_{i}$.
Then we have
$ \mbox{\textbf{\texttt{g}}}^{\mbox{\tiny \rm G}\!}(k)
=
  \mbox{\textbf{\texttt{g}}}_{0}^{\mbox{\tiny \rm G}\!}(k)
+ \mbox{\textbf{\texttt{g}}}_{1}^{\mbox{\tiny \rm G}\!}(k)
+ \mbox{\textbf{\texttt{g}}}_{2}^{\mbox{\tiny \rm G}\!}(k)
+ \mbox{\textbf{\texttt{g}}}_{3}^{\mbox{\tiny \rm G}\!}(k)$.
We can also obtain the equalities
$\mbox{\textbf{\texttt{g}}}_{0}^{\mbox{\tiny \rm G}\!}(k)
\!=\!
 \mbox{\textbf{\texttt{g}}}_{0}^{\mbox{\tiny \rm G}\!}(k\!-\!1)
\!+\!
 \mbox{\textbf{\texttt{g}}}_{1}^{\mbox{\tiny \rm G}\!}(k\!-\!1)$,
$\mbox{\textbf{\texttt{g}}}_{1}^{\mbox{\tiny \rm G}\!}(k)
\!=\!
 \mbox{\textbf{\texttt{g}}}_{3}^{\mbox{\tiny \rm G}\!}(k\!-\!1)$,
$\mbox{\textbf{\texttt{g}}}_{2}^{\mbox{\tiny \rm G}\!}(k)
\!=\!
 \mbox{\textbf{\texttt{g}}}_{2}^{\mbox{\tiny \rm G}\!}(k\!-\!1)$,
$\mbox{\textbf{\texttt{g}}}_{3}^{\mbox{\tiny \rm G}\!}(k)
\!=\!
 \mbox{\textbf{\texttt{g}}}_{0}^{\mbox{\tiny \rm G}\!}(k\!-\!1)
\!+\!
 \mbox{\textbf{\texttt{g}}}_{1}^{\mbox{\tiny \rm G}\!}(k\!-\!1)$.
Since the only child of $\nabla_{2}$ is $\nabla_{2}$, we can deduce
$\mbox{\textbf{\texttt{g}}}_{2}^{\mbox{\tiny \rm G}\!}(k)
\!=\!
 \mbox{\textbf{\texttt{g}}}_{2}^{\mbox{\tiny \rm G}\!}(k\!-\!1)
\!=\!
\cdots
\!=\!
\mbox{\textbf{\texttt{g}}}_{2}^{\mbox{\tiny \rm G}\!}(1)
\!=\! 1.$
Putting things together we obtain the recurrence for $\mbox{\textbf{\texttt{g}}}^{\mbox{\tiny \rm G}\!}(k)$:
$$
\begin{array}{rcl}
\mbox{\textbf{\texttt{g}}}^{\mbox{\tiny \rm G}\!}(k)
\!\!&\!\!=\!\!&\!\!
\mbox{\textbf{\texttt{g}}}_{0}^{\mbox{\tiny \rm G}\!}(k)
+
\mbox{\textbf{\texttt{g}}}_{1}^{\mbox{\tiny \rm G}\!}(k)
+
\mbox{\textbf{\texttt{g}}}_{2}^{\mbox{\tiny \rm G}\!}(k)
+
\mbox{\textbf{\texttt{g}}}_{3}^{\mbox{\tiny \rm G}\!}(k)                                             \\[1.25 ex]
\!\!&\!\!=\!\!&\!\!
\bigl(
\mbox{\textbf{\texttt{g}}}_{0}^{\mbox{\tiny \rm G}\!}(k\!-\!1)
\!+\!
\mbox{\textbf{\texttt{g}}}_{1}^{\mbox{\tiny \rm G}\!}(k\!-\!1)
\bigr)
\!+\!
\mbox{\textbf{\texttt{g}}}_{3}^{\mbox{\tiny \rm G}\!}(k\!-\!1)
\!+\!
\mbox{\textbf{\texttt{g}}}_{2}^{\mbox{\tiny \rm G}\!}(k\!-\!1)
\!+\!
\bigl(
\mbox{\textbf{\texttt{g}}}_{0}^{\mbox{\tiny \rm G}\!}(k\!-\!1)
\!+\!
\mbox{\textbf{\texttt{g}}}_{1}^{\mbox{\tiny \rm G}\!}(k\!-\!1)
\bigr)                                                                                               \\[1.25 ex]
\!\!&\!\!=\!\!&\!\!
\mbox{\textbf{\texttt{g}}}^{\mbox{\tiny \rm G}\!}(k\!-\!1)
\!+\!
\mbox{\textbf{\texttt{g}}}_{0}^{\mbox{\tiny \rm G}\!}(k\!-\!1)
\!+\!
\mbox{\textbf{\texttt{g}}}_{1}^{\mbox{\tiny \rm G}\!}(k\!-\!1)                                       \\[1.25 ex]
\!\!&\!\!=\!\!&\!\!
\mbox{\textbf{\texttt{g}}}^{\mbox{\tiny \rm G}\!}(k\!-\!1)
\!+\!
\bigl(
\mbox{\textbf{\texttt{g}}}_{0}^{\mbox{\tiny \rm G}\!}(k\!-\!2)
\!+\!
\mbox{\textbf{\texttt{g}}}_{1}^{\mbox{\tiny \rm G}\!}(k\!-\!2)
\bigr)
\!+\!
\mbox{\textbf{\texttt{g}}}_{3}^{\mbox{\tiny \rm G}\!}(k\!-\!2)
\!+\!
\mbox{\textbf{\texttt{g}}}_{2}^{\mbox{\tiny \rm G}\!}(k\!-\!2)
\!-\!
\mbox{\textbf{\texttt{g}}}_{2}^{\mbox{\tiny \rm G}\!}(k\!-\!2)                                       \\[1.25 ex]
\!\!&\!\!=\!\!&\!\!
\mbox{\textbf{\texttt{g}}}^{\mbox{\tiny \rm G}\!}(k\!-\!1)
\!+\!
\mbox{\textbf{\texttt{g}}}^{\mbox{\tiny \rm G}\!}(k\!-\!2) \!-\! 1.
\end{array}
$$
Substituting $\mbox{\textbf{\texttt{t}}}(k)\!=\!\mbox{\textbf{\texttt{g}}}^{\mbox{\tiny \rm G}\!}(k)-1$
into the previous formula we obtain the recurrence
$
\mbox{\textbf{\texttt{t}}$(k)$}\!=\!\mbox{\textbf{\texttt{t}}$(k\!-\!1)$}\break
+\mbox{\textbf{\texttt{t}}$(k\!-\!2)$}$.
Based on the initial conditions
$\mbox{\textbf{\texttt{g}}}^{\mbox{\tiny \rm G}\!}(1)\!=\!4$
and
$\mbox{\textbf{\texttt{g}}}^{\mbox{\tiny \rm G}\!}(2)\!=\!6$,
i.e.,
$\mbox{\mbox{\textbf{\texttt{t}}}(1)}\!=\!3$
and
$\mbox{\mbox{\textbf{\texttt{t}}}(2)}\!=\!5$,
we conclude that
$\mbox{\textbf{\texttt{g}}}^{\mbox{\tiny \rm G}\!}(k)
\!=\!
F_{k+3} \!+\! 1$.

\section{The enumerations of the higher order non-trivial compositions of the differential operations
         and the directional derivative in the space $\mathbb{R}^n$}

\smallskip
We start this section by recalling some definitions of the theory of differential forms.
Denote by $\mathbb{R}^n$ the $n$-dimensional Euclidean space ($n \!\geq\! 3$) and consider the set of smooth functions
$$
\mbox{\rm A}_{0} = \{ f\!:\! \mathbb{R}^n \! \rightarrow \! \mathbb{R} \, | \, f\!\in\!C^{\infty}(\mathbb{R}^n) \}.
$$
The set of all differential $k$-forms on $\mathbb{R}^n$, denoted by $\Omega^k(\mathbb{R}^n)$,
is a free $\mbox{\rm A}_{0}$-module of the rank ${n \choose k}$
with the standard basis $\{dx_I = dx_{i_1} \cdots dx_{i_k} \,|\, 1\!\leq\! i_1\!< \!\cdots \!<\! i_k \!\leq\! n\}$.
A differential $k$-form $\omega$ can be written uniquely as $\omega \!=\! \sum_{I \in \mathcal{I}} \omega_I dx_I$,
where $\omega_I \in \mbox{\rm A}_{0}$, and $\mathcal{I} \!=\! \mathcal{I}(k,n)$ is the set of multi-indices
$I \!=\! (i_1, \ldots, i_k)$, ($1 \!\leq\! i_1 \!<\! \cdots \!<\! i_k \!\leq\! n$).
The complement of $I$ is $J \!=\! (j_1, \ldots, j_{n-k}) \in \mathcal{I}(n-k,n)$, ($1\!\leq\! j_1\!<\! \cdots \!<\! j_{n-k} \!\leq\! n$),
where components $j_p$ are the elements of the set $\{1, \ldots, n\}\backslash \{i_1, \ldots, i_k\}$.
We have $dx_I dx_J \!=\! \sigma(I) dx_1 \ldots dx_n$, where $\sigma(I)$ is a signature of the permutation
$(i_1, \ldots, i_k$, $j_1, \ldots, j_{n-k})$. Note that $\sigma(J) \!=\! (-1)^{k(n-k)}\sigma(I)$.
With the notion mentioned above we define $\star_k (dx_I) \!=\! \sigma(I) dx_J$.
A map $\star_k\!:\!\Omega^k(\mathbb{R}^n) \!\longrightarrow\! \Omega^{n-k}(\mathbb{R}^n)$ defined by
$\star_k (\omega) \!=\! \sum_{I \in \mathcal{I}(k,n)} \omega_I \star_k (dx_I)$
is the Hodge star operator and it provides a natural isomorphism between $\Omega^k(\mathbb{R}^n)$
and $\Omega^{n-k}(\mathbb{R}^n)$. The Hodge star operator applied twice to a differential $k$-form yields
$\star_{n-k}(\star_k \omega) \!=\! (-1)^{nk+k+s}(\omega)$, where $s$ is the number of negative signs
in the inner product of the base vectors of the space $\mathbb{R}^n$ (see \cite{BalasubramanianLynnSenGupta}, p.$\,$29).
For the inverse of $\star_k$ the equality ${\star_k}^{\!-1}(\psi) \!=\! (-1)^{nk+k+s}\star_{n-k}(\psi)$ holds,
where \mbox{$\psi \!\in\! \Omega^{n-k}(\mathbb{R}^n)$}.

\smallskip

A differential $0$-form is a function $f(x_1, \ldots, x_n) \!\in\! \mbox{\rm A}_{0}$.
We define $df$ to be the differential $1$-form $df\!=\!\sum_{i=1}^n{\!\frac{\partial f}{\partial x_i}dx_i}$.
Given a differential $k$-form $\sum_{I \in \mathcal{I}} \omega_I dx_I$,
the exterior derivative $d_k \omega$ is the differential $(k\!+\!1)$-form
$d_k \omega = \sum_{I \in \mathcal{I}}d\omega_I dx_I$.
The exterior derivative $d_k$ is a linear map
$d_k\!:\!\Omega^k(\mathbb{R}^n) \!\longrightarrow\! \Omega^{k+1}(\mathbb{R}^n)$
which obeys Leibnitz rule
$$d_{p+q}(\omega \, \psi) \!=\! (d_p \omega) \, \psi + (-1)^p \omega \, (d_q\psi),$$
$\omega$ and $\psi$ are differential $p\,$-form and $q\,$-form.
The exterior derivative has a property
$$d_{k+1}(d_k \omega) = 0,$$
for any differential $k$-form $\omega$.
For more details on the topic please refer to \cite{BalasubramanianLynnSenGupta, PerotZusi14, Weintraub96}.
Some historical notes on the development of the theory of differential forms are given in \cite{Katz85, Walter07}.

\smallskip

Consider the sets of the smooth functions
$$
\mbox{\rm A}_{k}
\!=\!
\mbox{\large $\{$}
\vec{f}\!:\! \mathbb{R}^n \!\rightarrow\! \mathbb{R}^{n \choose k} \, | \, f_1, \ldots, f_{{n \choose k}} \!\in\!C^{\infty}(\mathbb{R}^n)
\mbox{\large $\}$}
\quad
\bigl(
\,m =\! \left\lfloor \mbox{\small $\displaystyle\frac{n}{2}$} \right\rfloor\!,\,0 \leq k \leq m\,
\bigr).
$$

\break

\noindent
Let $p_k \!:\! \Omega^{k}(\mathbb{R}^{n}) \!\rightarrow\! \mbox{\rm A}_{k}$
be a presentation of the differential form in coordinate notation.
Define functions
$\varphi_{i}$ $(0 \!\leq\! i \!\leq\! m)$ and $\varphi_{n-j}$ $(0 \leq j < n\!-\!m)$ as follows:
$$
\varphi_{i}
=
p_{i} \!:\! \Omega^{i}(\mathbb{R}^{n}) \longrightarrow \mbox{\rm A}_{i}
\quad
\mbox{and}
\quad
\varphi_{n-j}
=
p_{j} \, \star_{j}^{-1} \!:\! \Omega^{n-j}(\mathbb{R}^{n}) \longrightarrow \mbox{\rm A}_{j}.
$$
Analogously to Male\v sevi\' c \cite{HiOrd06}, the combination of the Hodge star operator and the exterior derivative
generates differential operations:
$$
\nabla_{k} = \varphi_{k} \, d_{k-1} \, \varphi_{k-1}^{-1} \quad (1\!\leq\!k\!\leq\! n).
$$
Please see the list below.
$$
\begin{array}{c}
\begin{tabular}{cc}
$\begin{array}{ll}
\hspace*{-7.5 mm}   \mbox{\small $\mbox{$\mathcal{A}$}_{n}\;(n\!=\!2\,m)\,$:}   \hspace*{-60.0 mm}                                                  \\[ 0.300 ex]
\hspace*{-7.0 mm} & \mbox{\footnotesize$\nabla_{1}\!=\!p_{1}\,d_{0}\,p_{0}^{-1}\!:\!A_{0}\!\rightarrow\!A_{1}$}                                     \\[ 0.000 ex]
\hspace*{-7.0 mm} & \mbox{\footnotesize$\nabla_{2}\!=\!p_{2}\,d_{1}\,p_{1}^{-1}\!:\!A_{1}\!\rightarrow\!A_{2}$}                                     \\[-0.875 ex]
\hspace*{-7.0 mm} & \,\,\mbox{\small $\vdots$}                                                                                                      \\[-0.275 ex]
\hspace*{-7.0 mm} & \mbox{\footnotesize$\nabla_{i}\!=\!p_{i}\,d_{i-1}\,p_{i-1}^{-1}\!:\!A_{i-1}\!\rightarrow\!A_{i}$}                               \\[-0.875 ex]
\hspace*{-7.0 mm} & \,\,\mbox{\small $\vdots$}                                                                                                      \\[-0.275 ex]
\hspace*{-7.0 mm} & \mbox{\footnotesize$\nabla_{m}\!=\!p_{m}\,d_{m-1}\,p_{m-1}^{-1}\!:\!A_{m-1}\!\rightarrow\!A_{m}$}                               \\[ 0.000 ex]
\hspace*{-7.0 mm} & \mbox{\footnotesize$\nabla_{m+1}\!=\!p_{m-1}\star_{m-1}^{-1}d_{m}\,p_{m}^{-1}\!:\!A_{m}\!\rightarrow\!A_{m-1}$}                 \\[ 0.000 ex]
\hspace*{-7.0 mm} & \mbox{\footnotesize$\nabla_{m+2}\!=\!p_{m-2}\star_{m-2}^{-1}d_{m+1}\star_{m-1}p_{m-1}^{-1}\!:\!A_{m-1}\!\rightarrow\!A_{m-2}$}  \\[-0.875 ex]
\hspace*{-7.0 mm} & \,\,\mbox{\small $\vdots$}                                                                                                      \\[-0.275 ex]
\hspace*{-7.0 mm} & \mbox{\footnotesize$\nabla_{n-j}\!=\!p_{j}\star_{j}^{-1}d_{n-(j+1)}\star_{j+1}p_{j+1}^{-1}\!:\!A_{j+1}\!\rightarrow\!A_{j}$}    \\[-0.875 ex]
\hspace*{-7.0 mm} & \,\,\mbox{\small $\vdots$}                                                                                                      \\[-0.275 ex]
\hspace*{-7.0 mm} & \mbox{\footnotesize$\nabla_{n-1}\!=\!p_{1}\star_{1}^{-1}d_{n-2}\star_{2}p_{2}^{-1}\!:\!A_{2}\!\rightarrow\!A_{1}$}              \\[ 0.000 ex]
\hspace*{-7.0 mm} & \mbox{\footnotesize$\nabla_{n}\!=\!p_{0}\star_{0}^{-1}d_{n-1}\star_{1}p_{1}^{-1}\!:\!A_{1}\!\rightarrow\!A_{0}$}
\,\mbox{\normalsize ,}
\end{array}$
&
$\hspace*{8.0 mm}
\begin{array}{ll}
\hspace*{-15.0 mm}   \mbox{\small $\mbox{$\mathcal{A}$}_{n}\;(n\!=\!2\,m\!+\!1)\,$:}    \hspace*{-60.0 mm}                                          \\[ 0.300 ex]
\hspace*{-14.0 mm} & \mbox{\footnotesize$\nabla_{1}\!=\!p_{1}\,d_{0}\,p_{0}^{-1}\!:\!A_{0}\!\rightarrow\!A_{1}$}                                    \\[ 0.000 ex]
\hspace*{-14.0 mm} & \mbox{\footnotesize$\nabla_{2}\!=\!p_{2}\,d_{1}\,p_{1}^{-1}\!:\!A_{1}\!\rightarrow\!A_{2}$}                                    \\[-0.875 ex]
\hspace*{-14.0 mm} & \,\,\mbox{\small $\vdots$}                                                                                                     \\[-0.275 ex]
\hspace*{-14.0 mm} & \mbox{\footnotesize$\nabla_{i}\!=\!p_{i}\,d_{i-1}\,p_{i-1}^{-1}\!:\!A_{i-1}\!\rightarrow\!A_{i}$}                              \\[-0.875 ex]
\hspace*{-14.0 mm} & \,\,\mbox{\small $\vdots$}                                                                                                     \\[-0.275 ex]
\hspace*{-14.0 mm} & \mbox{\footnotesize$\nabla_{m}\!=\!p_{m}\,d_{m-1}\,p_{m-1}^{-1}\!:\!A_{m-1}\!\rightarrow\!A_{m}$}                              \\[ 0.000 ex]
\hspace*{-14.0 mm} & \mbox{\footnotesize$\nabla_{m+1}\!=\!p_{m}\star_{m}^{-1}d_{m}\,p_{m}^{-1}\!:\!A_{m}\!\rightarrow\!A_{m}$}                      \\[ 0.000 ex]
\hspace*{-14.0 mm} & \mbox{\footnotesize$\nabla_{m+2}\!=\!p_{m-1}\star_{m-1}^{-1}d_{m+1}\star_{m}p_{m}^{-1}\!:\!A_{m}\!\rightarrow\!A_{m-1}$}       \\[ 0.000 ex]
\hspace*{-14.0 mm} & \mbox{\footnotesize$\nabla_{m+3}\!=\!p_{m-2}\star_{m-2}^{-1}d_{m+2}\star_{m-1}p_{m-1}^{-1}\!:\!A_{m-1}\!\rightarrow\!A_{m-2}$} \\[-0.875 ex]
\hspace*{-14.0 mm} & \,\,\mbox{\small $\vdots$}                                                                                                     \\[-0.275 ex]
\hspace*{-14.0 mm} & \mbox{\footnotesize$\nabla_{n-j}\!=\!p_{j}\star_{j}^{-1}d_{n-(j+1)}\star_{j+1}p_{j+1}^{-1}\!:\!A_{j+1}\!\rightarrow\!A_{j}$}   \\[-0.875 ex]
\hspace*{-14.0 mm} & \,\,\mbox{\small $\vdots$}                                                                                                     \\[-0.275 ex]
\hspace*{-14.0 mm} & \mbox{\footnotesize$\nabla_{n-1}\!=\!p_{1}\star_{1}^{-1}d_{n-2}\star_{2}p_{2}^{-1}\!:\!A_{2}\!\rightarrow\!A_{1}$}             \\[ 0.000 ex]
\hspace*{-14.0 mm} & \mbox{\footnotesize$\nabla_{n}\!=\!p_{0}\star_{0}^{-1}d_{n-1}\star_{1}p_{1}^{-1}\!:\!A_{1}\!\rightarrow\!A_{0}$}
\,\mbox{\normalsize .}
\end{array}$
\end{tabular}                                                                                                                                       \\
\\[-2.00 ex]
\hspace*{-16.00 mm}
\mbox{\small Table 1. List of differential operations in $\mathbb{R}^n$}
\end{array}
$$
For $n\!=\!3$, we obtain the standard definitions of the gradient, curl and divergence,
see \cite{BalasubramanianLynnSenGupta, Kotiuga89}. For $n\!=\!4$ let us consider Lewis--Wilson
four-dimensional differential operations
\cite{ENWilsonGNLewis1912}:
$$
\begin{array}{rcl}
\nabla_1(f)
\!&\!\!=\!\!&\!
\left(
\mbox{\footnotesize $\displaystyle\frac{\partial f}{\partial x_1}$},
\mbox{\footnotesize $\displaystyle\frac{\partial f}{\partial x_2}$},
\mbox{\footnotesize $\displaystyle\frac{\partial f}{\partial x_3}$},
\mbox{\small $-$}
\mbox{\footnotesize $\displaystyle\frac{\partial f}{\partial x_4}$}
\right)
:
A_{0} \longrightarrow A_{1},                                                   \\[1.5 ex]
\nabla_2{\big (}\mbox{\bf f}{\big )}
\!&\!\!=\!\!&\!
\nabla_2{\big (}(\mbox{\small $f_1$},\mbox{\small $f_2$},\mbox{\small $f_3$},\mbox{\small $f_4$}){\big )}
                                                                               \\[0.25 ex]
\!&\!\!=\!\!&\!
{\Big (}
\mbox{\footnotesize $\displaystyle\frac{\partial f_3}{\partial x_2}$}
\!-\!
\mbox{\footnotesize $\displaystyle\frac{\partial f_2}{\partial x_3}$},
\mbox{\footnotesize $\displaystyle\frac{\partial f_1}{\partial x_3}$}
\!-\!
\mbox{\footnotesize $\displaystyle\frac{\partial f_3}{\partial x_1}$},
\mbox{\footnotesize $\displaystyle\frac{\partial f_2}{\partial x_1}$}
\!-\!
\mbox{\footnotesize $\displaystyle\frac{\partial f_1}{\partial x_2}$},
\mbox{\footnotesize $\displaystyle\frac{\partial f_4}{\partial x_1}$}
\!+\!
\mbox{\footnotesize $\displaystyle\frac{\partial f_1}{\partial x_4}$},
\mbox{\footnotesize $\displaystyle\frac{\partial f_4}{\partial x_2}$}
\!+\!
\mbox{\footnotesize $\displaystyle\frac{\partial f_2}{\partial x_4}$},
\mbox{\footnotesize $\displaystyle\frac{\partial f_4}{\partial x_3}$}
\!+\!
\mbox{\footnotesize $\displaystyle\frac{\partial f_3}{\partial x_4}$}
{\Big )}
:
A_{1} \longrightarrow A_{2},                                                   \\[1.5 ex]
\nabla_3{\big (}\mbox{\bf F}{\big )}
\!\!&\!\!=\!\!&\!
\nabla_3{\big (}(\mbox{\small $F_1$},\mbox{\small $F_2$},\mbox{\small $F_3$},\mbox{\small $F_4$},\mbox{\small $F_5$},\mbox{\small $F_6$}){\big )}
                                                                               \\[0.25 ex]
\!&\!=\!&\!
{\Big (}
\mbox{\footnotesize $\displaystyle\frac{\partial F_6}{\partial x_2}$}
\!-\!
\mbox{\footnotesize $\displaystyle\frac{\partial F_5}{\partial x_3}$}
\!-\!
\mbox{\footnotesize $\displaystyle\frac{\partial F_1}{\partial x_4}$},
\mbox{\footnotesize $\displaystyle\frac{\partial F_4}{\partial x_3}$}
\!-\!
\mbox{\footnotesize $\displaystyle\frac{\partial F_6}{\partial x_1}$}
\!-\!
\mbox{\footnotesize $\displaystyle\frac{\partial F_2}{\partial x_4}$},
\mbox{\footnotesize $\displaystyle\frac{\partial F_5}{\partial x_1}$}
\!-\!
\mbox{\footnotesize $\displaystyle\frac{\partial F_4}{\partial x_2}$}
\!-\!
\mbox{\footnotesize $\displaystyle\frac{\partial F_3}{\partial x_4}$},
\mbox{\footnotesize $\displaystyle\frac{\partial F_1}{\partial x_1}$}
\!+\!
\mbox{\footnotesize $\displaystyle\frac{\partial F_2}{\partial x_2}$}
\!+\!
\mbox{\footnotesize $\displaystyle\frac{\partial F_3}{\partial x_3}$}
{\Big )}
:
A_{2} \longrightarrow A_{1},                                                   \\[1.5 ex]
\nabla_4{\big (}\mbox{\bf f}{\big )}
\!\!&\!\!=\!\!&\!
\nabla_4{\big (}(\mbox{\small $f_1$},\mbox{\small $f_2$},\mbox{\small $f_3$},\mbox{\small $f_4$}){\big )}
\;=\;
\mbox{\footnotesize $\displaystyle\frac{\partial f_1}{\partial x_1}$}
\!+\!
\mbox{\footnotesize $\displaystyle\frac{\partial f_2}{\partial x_2}$}
\!+\!
\mbox{\footnotesize $\displaystyle\frac{\partial f_3}{\partial x_3}$}
\!+\!
\mbox{\footnotesize $\displaystyle\frac{\partial f_4}{\partial x_4}$}
:
A_{1} \longrightarrow A_{0},
\end{array}
$$
over the set $\mathcal{A}_{4} =\! \{\nabla_1, \nabla_2, \nabla_3, \nabla_4 \}$. Using Lewis--Wilson determination
of the Hodge operators on the base vectors (see \cite{ENWilsonGNLewis1912}, p.$\,$450) we may conclude that Lewis--Wilson
determination of four-dimensional differential operations is consentient with previous definition of the differential
operations in the space $\mathbb{R}^4$. Therefore $\nabla_2 \circ \nabla_1 \!=\! 0$ and $\nabla_3 \circ \nabla_2 \!=\! 0$
and $\nabla_4 \circ \nabla_3 \!=\! 0$~are~true. Let~us~emphasize that, analogously to results of F.C. Chang~\cite{FC_Chang05, FC_Chang12},
previous four-dimensional differential operations could be determined by the appropriate matrix formulation.

Let us note that the formulas for determination a number of the compositions over the set $\mathcal{A}_{n}$ and
the corresponding recurrences are obtained by Male\v sevi\' c \cite{HiOrd98, HiOrd06}.
An application of the formulas is given by Myers \cite{Myers}.
The corresponding integer sequences can be found in \cite{Sloane}.
For the proofs of the following two theorems we refer the reader to \cite{HiOrd98}.

\begin{theorem}
\label{non trivial}
A non-trivial composition over the set $\mathcal{A}_{n}$ is of the form:
$$
( \, \nabla_i \, \circ \, ) \, \nabla_{n+1-i} \circ \nabla_i \circ \, \cdots \, \circ \nabla_{n+1-i} \circ \nabla_i,
$$
for some $i$ $(2 i, \, 2i\!-\!2 \neq n, \; 1 \leq i \leq n)$.
The term in brackets is included if the number of the differential operations is odd and is left out otherwise.
\end{theorem}

In the terms of the Hodge star operator and the exterior derivative, we have the following representation of the non-trivial composition
$$
\begin{array}{l}
(\, \nabla_i \, \circ \,) \, \nabla_{n+1-i} \circ \nabla_i \circ \, \cdots \, \circ \nabla_{n+1-i} \circ \nabla_i  =    \\[2.0 ex]
\begin{cases}
\mbox{\small $\!(\, p_i \, d_{i-1} \,  p_{i-1}^{-1} \, ) \, p_{i-1} \star_{i-1}^{-1} \, d_{n-i} \, \star_i \, d_i
\, \cdots \, \star_{i-1}^{-1} \, d_{n-i} \, \star_i \, d_{i-1} \; p_{i-1}$},
&  \mbox{\small $i \leq m$};                                                                                            \\[1.5 ex]
\!\!\!\begin{array}{l}
\mbox{\small $\!(\, p_{n-i} \, \star_{n-i}^{-1} \, d_{i-1} \, \star_{n+1-i} \, p_{n+1-i}^{-1} \,)
\, p_{n+1-i} \, d_{n-i} \, \star_{n-i}^{-1} \, d_{i-1} \, \star_{n+1-i}$}                                               \\[0.5 ex]
\mbox{\small $\,\cdots \, d_{n-i} \, \star_{n-i}^{-1} \, d_{i-1} \, \star_{n+1-i} \, p_{n+1-i}^{-1}$},
\end{array}
&  \mbox{\small $i > m$};
\end{cases}
\end{array}
$$
where $i \!\in\! \{1, \ldots, n\} \backslash \{m, m+1\}$ if $n \!=\! 2\,m\,$ or $i \!\in\! \{1, \ldots, n\} $ if $n \!=\! 2\,m + 1$.
\begin{theorem}
\label{non-trivial A_n}
Let $\mbox{\textbf{\tt g}}(k)$ be the number of the $k^{\mbox{\scriptsize \it th}}$ order
non-trivial compositions over the set $\mathcal{A}_{n}$.$\;$Then we have
$$
\mbox{\textbf{\tt g}}(k)
=
\begin{cases}
n,   &      \!         2 \nmid n \; ; \\[0.50 ex]
n,   & 2 \, | \, n \:\:,\:\: k=1 \; ; \\[0.50 ex]
n-1, & 2 \, | \, n \:\:,\:\: k=2 \; ; \\[0.50 ex]
n-2, & 2 \, | \, n \:\:,\:\: k>2 \; .
\end{cases}
$$
\end{theorem}

The Hodge dual to the exterior derivative $d_k \!:\! \Omega^k(\mathbb{R}^n) \rightarrow \Omega^{k+1}(\mathbb{R}^n)$
is codifferential $\delta_{k-1}$. It is a generalization of the divergence.
The codifferential is a linear map $\delta_{k-1} \!:\! \Omega^k(\mathbb{R}^n) \rightarrow \Omega^{k-1}(\mathbb{R}^n)$,
determined by
$$
\delta_{k-1} = (-1)^{nk+k+s+1} \star_{n-(k-1)} d_{n-k} \star_{k} = (-1)^k \star_{k-1}^{-1} d_{n-k} \star_{k}
$$
(see \cite{BalasubramanianLynnSenGupta}, p.$\,$33). Note that $\nabla_{n-j} = (-1)^{j+1}\,p_{j}\;\delta_j\;p_{j+1}^{-1}$,
($0 \leq j < n \!-\! m \!-\! 1$). The codifferential~can~be coupled with the exterior derivative to construct
the Hodge Laplacian, also known as the Laplace-de Rham operator,
$\Delta_k \!:\! \Omega^k(\mathbb{R}^n) \rightarrow \Omega^k(\mathbb{R}^n)$.
The Hodge Laplacian is a harmonic generalization of the Laplace differential operator, given by
$\Delta_0 = \delta_0 d_0$ and \mbox{$\Delta_k = \delta_k d_k \!+\!  d_{k-1}\delta_{k-1}$},
for $1 \!\leq\! k \!\leq\! m$, see \cite{vonWestenholz78}.
The operator $\Delta_0$ is actually the negative of the Laplace-Beltrami (scalar) operator.

\smallskip

A $k$-form $\omega$ is called harmonic if $\Delta_k (\omega) = 0$.
We say that $\vec{f} \in A_{k}$ is a harmonic function if $\omega = {p_k}^{-1}(\vec{f})$ is a harmonic $k$-form.
If $k \geq 1$ harmonic function $\vec{f}$ is also called a harmonic field.
For the function $\vec{f} \!\in\! A_{k}$ ($1 \!\leq\! k \!\leq\! m$) we have that
$\Delta_{k}(p_{k}^{-1}\!\vec{f})   \!=\! 0$
iff
$\delta_{k-1}(p_{k}^{-1}\!\vec{f}) \!=\! 0$
and
$d_{k}(p_{k}^{-1}\!\vec{f})        \!=\! 0$,
\cite[Proposition 4.15]{vonWestenholz78}.
In fact, we obtain the following lemma.

\begin{lemma}
Let $\vec{f} \in A_{k}$ $(1 \!\leq\! k \!\leq\! m)$.$\;$Then
$$
\Delta_{k}(p_{k}^{-1}\!\vec{f}) = 0
\;\Longleftrightarrow\;
\nabla_{n-(k-1)}(\vec{f}) = 0
\,\wedge\,
\nabla_{k+1}(\vec{f}) = 0.
$$
\end{lemma}

For harmonic function $f \!\in\! A_{0}$ we have $\Delta_0 f \!=\! \delta_0 d_0 f \!=\! 0$.
Hence we get \mbox{$\nabla_n \!\circ\! \nabla_1 f \!=\! 0$}
and, consequently, we obtain $(\nabla_1 \circ) \nabla_n \circ \nabla_1 \circ \dots \circ \nabla_n \circ \nabla_1 f \!=\! 0$.
We can now rephrase Theorem \ref{non trivial} for harmonic functions.

\begin{theorem}
All the second and higher order non-trivial compositions over the set $\mathcal{A}_{n}$
acting on harmonic function $f \!\in\! A_{0}$ are trivial. Furthermore,
all the first and higher order non-trivial compositions over the set $\mathcal{A}_{n}$
acting on harmonic field $\vec{f} \!\in\! A_{k}$ $(1 \!\leq\! k \!\leq\! m)$ are trivial.
\end{theorem}

We say that $\vec{f} \!\in\! A_{k}$ ($1 \!\leq\! k \!\leq\! m$) is a coordinate-harmonic function or that
$\vec{f}$ satisfies harmonic coordinate condition, if all its coordinates are harmonic functions.
Male\v sevi\' c \cite{HiOrd96} showed that all the third and higher order non-trivial compositions of the differential operations
acting on coordinate-harmonic functions are trivial in $\mathbb{R}^{3}$.

\begin{conjecture}
\label{Conjecture}
All the third and higher order non-trivial compositions over the set $\mathcal{A}_{n}$
acting on coordinate-harmonic functions are trivial in $\mathbb{R}^{n}$.
\end{conjecture}

Lewis and Wilson \cite{GNLewis1910, ENWilsonGNLewis1912} gave an approach to a coordinate investigation
of Conjecture \ref{Conjecture} for $n=4$, see also \cite{Sommerfeld10a, Sommerfeld10b, Kraft1911}.
A similar problem for coordinate-harmonic functions can be found in discrete exterior calculus \cite{AADM}
and combinatorial Hodge theory \cite{BAMS}. Remark that some applications of directional derivatives in discrete approximations
of higher order differential operations are considered in \cite{BelyaevKhesinTabachnikov}.

Let $f \!\in\! \mbox{\rm A}_{0}$ be a scalar function, and $\vec{e} = (e_1,\dots,e_n) \in \mathbb{R}^{n}$ be a unit vector.
The Gateaux directional derivative in a direction $\vec{e}$ is defined by
$$
\mbox{\dire} \, f = \nabla_0 f =
\displaystyle \sum\limits_{k=1}^{n}{ \frac{\partial f}{\partial x_k} \, e_k} \!:\! A_{0} \longrightarrow A_{0} \, .
$$
Extend a set of differential operations $\mathcal{A}_{n} = \{ \nabla_{1}, \dots, \nabla_{n} \}$
with the directional derivative $\nabla_0$ to the set
$\mathcal{B}_{n} \!=\! \mathcal{A}_{n} \!\cup \{ \nabla_{0} \} = \{ \nabla_{0}, \nabla_{1}, \dots, \nabla_{n} \}$.
Male\v sevi\' c, and Jovovi\' c \cite{HiOrd07} have given the recurrences for counting the number of compositions over the set $\mathcal{B}_{n}$.
The corresponding integer sequences can be find in \cite{Sloane}.

\smallskip

Consider the non-trivial compositions over the set $\mathcal{B}_n$ containing $\nabla_0$.
First of all, we can compose the directional derivative $\nabla_0$ by itself to obtain the non-trivial compositions.
In other words, the non-trivial compositions containing just $\nabla_0$ are obtained from $\nabla_0$
by substituting $\nabla_0 \mapsto \nabla_0 \circ \nabla_0$. We define the coordinate  Gateaux directional~derivative~by
\mbox{${\bm \nabla}_0(\vec{f}) \,\mbox{\small $\bm =$}\,
(\nabla_0(f_1), \ldots , \nabla_0(f_n)) \!:\! A_1 \longrightarrow A_1$},
for the unit vector $\,\vec{e}\,$. The~following~two~equalities
$$
\nabla_1 \circ \nabla_0 \,\mbox{\small $\bm =$}\, {\bm \nabla}_0 \circ \nabla_1
\qquad\mbox{and}\qquad
\nabla_0 \circ \nabla_n \,\mbox{\small $\bm =$}\, \nabla_n \circ {\bm \nabla}_0
$$
hold. Therefore, the non-trivial compositions over the set $\mathcal{B}_n$ containing $\nabla_0$ can be obtained from
$( \, \nabla_1 \circ \, ) \, \nabla_{n} \circ \nabla_1 \circ \, \cdots \, \circ \nabla_{n} \circ \nabla_1( \, \circ \nabla_n \,)$
by substituting $\nabla_1 \mapsto \nabla_1 \circ \nabla_0$ or $\nabla_n \mapsto \nabla_0 \circ \nabla_n$.
Summarizing these facts, we have the following theorem.
\begin{theorem}
A non-trivial composition over the set $\mathcal{B}_{n}$ has the one of the following forms:

\bigskip
\noindent
{\boldmath $\;\;\;(i)\;\;\,$}
$\nabla_0^k=\underbrace{\nabla_0 \circ \nabla_0 \circ \, \cdots \circ \, \nabla_0}_k \;\;\; (k  \in \mathbb{N})$;          \\[-0.5 ex]

\noindent
{\boldmath $\;\;(ii)\;\;$}
$( \, \nabla_i \, \circ \, ) \, \nabla_{n+1-i} \circ \nabla_i \circ \, \cdots \, \circ \nabla_{n+1-i} \circ \nabla_i \;\;\;
         ( 2 i, \, 2i\!-\!2 \neq n, \; 1 \leq i \leq n)$;                                                                   \\[0.5 ex]

\noindent
{\boldmath $\;\,(iii)\;\,$}
$( \, \nabla_1 \circ \, ) \, \nabla_0^{k_p} \circ \nabla_{n} \circ \nabla_1 \circ \, \cdots \,\circ \nabla_0^{k_2}
  \circ \nabla_{n} \circ \nabla_1 \circ \nabla_0^{k_1}  \;\;\; (k_1, \ldots , k_p \in \mathbb{N}\cup\{0\})$;              \\[0.5 ex]

\noindent
{\boldmath $\;\;(iv)\;\,$}
$( \, \nabla_0^{k_q} \circ \nabla_{n} \circ \, ) \, \nabla_1 \circ \nabla_0^{k_{q-1}} \circ \nabla_{n} \circ \, \cdots
          \, \circ \nabla_1 \circ \nabla_0^{k_1} \circ \nabla_n \;\;\; (k_1, \ldots , k_q \in \mathbb{N}\cup\{0\})$.

\smallskip

\bigskip
\noindent
The terms in brackets are included if the number of the differential operations is odd and is left out otherwise.
\end{theorem}
The number of the higher order non-trivial compositions over the set $\mathcal{ B}_{n}$
is determined by the binary relation $\nu$, defined by:
$$
\nabla_{\!i}  \nu  \nabla_{\!j} \; \mbox{iff} \;
(i\!=\!0 \wedge j\!=\!0) \vee (i\!=\!0 \wedge j\!=\!1) \vee
(i\!=\!n \wedge j\!=\!0) \vee (i\!+\!j\!=\!n\!+\!1 \wedge 2i \not = n).
$$
By applying Theorem \ref{non-trivial A_n}
we conclude that the number of the $k^{\mbox{\scriptsize \rm th}}$ order non-trivial compositions
starting with $\nabla_{2}$,\dots, $\nabla_{n-1}$ can be expressed by formula
$$
\mbox{\textbf{\texttt{j}}}(k)
=
\mbox{\textbf{\texttt{g}}}(k)-2
=
\begin{cases}
n-2, &        \!       2 \nmid n \; ; \\[0.50 ex]
n-2, & 2 \, | \, n \:\:,\:\: k=1 \; ; \\[0.50 ex]
n-3, & 2 \, | \, n \:\:,\:\: k=2 \; ; \\[0.50 ex]
n-4, & 2 \, | \, n \:\:,\:\: k>2 \; .
\end{cases}
$$
Let $\mbox{\textbf{\texttt{g}}}^{\mbox{\tiny \rm G}\!}(k)$ be the number of
the $k^{\mbox{\scriptsize \rm th}}$ order non-trivial compositions over the set $\mathcal{B}_{n}$.
Let $\mbox{\textbf{\texttt{g}}}_{0}^{\mbox{\tiny \rm G}\!}(k)$, $\mbox{\textbf{\texttt{g}}}_{1}^{\mbox{\tiny \rm G}\!}(k)$ and
$\mbox{\textbf{\texttt{g}}}_{n}^{\mbox{\tiny \rm G}\!}(k)$ be the numbers of the $k^{\mbox{\scriptsize \rm th}}$ order non-trivial
compositions starting with $\nabla_{0}$, $\nabla_{1}$ and $\nabla_{n}$, respectively. Then we have
$\mbox{\textbf{\texttt{g}}}^{\mbox{\tiny \rm G}\!}(k)
=
\mbox{\textbf{\texttt{g}}}_{0}^{\mbox{\tiny \rm G}\!}(k)
+
\mbox{\textbf{\texttt{g}}}_{1}^{\mbox{\tiny \rm G}\!}(k)
+
\mbox{\textbf{\texttt{j}}}(k)
+
\mbox{\textbf{\texttt{g}}}_{n}^{\mbox{\tiny \rm G}\!}(k).$
Denote
\mbox{$\widetilde{\mbox{\textbf{\texttt{g}}}}^{\mbox{\tiny \rm G}\!}(k)
= \mbox{\textbf{\texttt{g}}}_{0}^{\mbox{\tiny \rm G}\!}(k)
+ \mbox{\textbf{\texttt{g}}}_{1}^{\mbox{\tiny \rm G}\!}(k)
+ \mbox{\textbf{\texttt{g}}}_{n}^{\mbox{\tiny \rm G}\!}(k)$}.
The following three recurrences are true
$\mbox{\textbf{\texttt{g}}}_{0}^{\mbox{\tiny \rm G}\!}(k)
=
\mbox{\textbf{\texttt{g}}}_{0}^{\mbox{\tiny \rm G}\!}(k\!-\!1)
+
\mbox{\textbf{\texttt{g}}}_{1}^{\mbox{\tiny \rm G}\!}(k\!-\!1)$,
$\mbox{\textbf{\texttt{g}}}_{1}^{\mbox{\tiny \rm G}\!}(k)
=
\mbox{\textbf{\texttt{g}}}_{n}^{\mbox{\tiny \rm G}\!}(k\!-\!1)$,
$\mbox{\textbf{\texttt{g}}}_{n}^{\mbox{\tiny \rm G}\!}(k)
=
\mbox{\textbf{\texttt{g}}}_{0}^{\mbox{\tiny \rm G}\!}(k\!-\!1)
+
\mbox{\textbf{\texttt{g}}}_{1}^{\mbox{\tiny \rm G}\!}(k\!-\!1)$.
Thus, the recurrence for
$\widetilde{\mbox{\textbf{\texttt{g}}}}^{\mbox{\tiny \rm G}\!}(k)$ is of the form:
$$
\begin{array}{rclclcl}
\widetilde{\mbox{\textbf{\texttt{g}}}}^{\mbox{\tiny \rm G}\!}(k)
    \!&\!\!=\!\!&\!
\mbox{\textbf{\texttt{g}}}_{0}^{\mbox{\tiny \rm G}\!}(k)
+
\mbox{\textbf{\texttt{g}}}_{1}^{\mbox{\tiny \rm G}\!}(k)
+
\mbox{\textbf{\texttt{g}}}_{n}^{\mbox{\tiny \rm G}\!}(k)                                            \\[1.25 ex]
\!&\!\!=\!\!&\!
{\big (}
\mbox{\textbf{\texttt{g}}}_{0}^{\mbox{\tiny \rm G}\!}(k\!-\!1)
+
\mbox{\textbf{\texttt{g}}}_{1}^{\mbox{\tiny \rm G}\!}(k\!-\!1){\big )}
+
\mbox{\textbf{\texttt{g}}}_{n}^{\mbox{\tiny \rm G}\!}(k\!-\!1)
+
{\big (}\mbox{\textbf{\texttt{g}}}_{0}^{\mbox{\tiny \rm G}\!}(k\!-\!1)
+
\mbox{\textbf{\texttt{g}}}_{1}^{\mbox{\tiny \rm G}\!}(k\!-\!1){\big )}                              \\[1.25 ex]
\!&\!\!=\!\!&
\widetilde{\mbox{\textbf{\texttt{g}}}}^{\mbox{\tiny \rm G}\!}(k\!-\!1)
+
\mbox{\textbf{\texttt{g}}}_{0}^{\mbox{\tiny \rm G}\!}(k\!-\!1)
+
\mbox{\textbf{\texttt{g}}}_{1}^{\mbox{\tiny \rm G}\!}(k\!-\!1)                                      \\[1.25 ex]
\!&\!\!=\!\!&\!
\widetilde{\mbox{\textbf{\texttt{g}}}}^{\mbox{\tiny \rm G}\!}(k\!-\!1)
+
{\big (}\mbox{\textbf{\texttt{g}}}_{0}^{\mbox{\tiny \rm G}\!}(k\!-\!2)
+
\mbox{\textbf{\texttt{g}}}_{1}^{\mbox{\tiny \rm G}\!}(k\!-\!2){\big )}
+
\mbox{\textbf{\texttt{g}}}_{n}^{\mbox{\tiny \rm G}\!}(k\!-\!2)                                      \\[1.25 ex]
\!&\!\!=\!\!&\!
\widetilde{\mbox{\textbf{\texttt{g}}}}^{\mbox{\tiny \rm G}\!}(k\!-\!1)
+
\widetilde{\mbox{\textbf{\texttt{g}}}}^{\mbox{\tiny \rm G}\!}(k\!-\!2).
\end{array}
$$
With initial conditions
$\widetilde{\mbox{\textbf{\texttt{g}}}}^{\mbox{\tiny \rm G}\!}(1)\!=\!3$,
$\widetilde{\mbox{\textbf{\texttt{g}}}}^{\mbox{\tiny \rm G}\!}(2)\!=\!5$
we deduce
$\widetilde{\mbox{\textbf{\texttt{g}}}}^{\mbox{\tiny \rm G}\!}(k)\!=\!F_{k+3}$.
Therefore, we have proved the following theorem.

\begin{theorem}
The number of the $k^{\mbox{\scriptsize \it th}}$ order non-trivial compositions over the set $\mathcal{B}_{n}$ is

\smallskip
\noindent
$$
\mbox{\textbf{\tt g}}^{\mbox{\tiny \rm G}\!}(k)
=
F_{k+3} + \mbox{\textbf{\tt j}}(k)
=
\begin{cases}
F_{k+3} + n-2 ,&                   \!   2 \nmid n \; ; \\[1.0 ex]
n+1           ,&      2   \, | \, n \:\:,\:\: k=1 \; ; \\[1.0 ex]
n+2           ,&      2   \, | \, n \:\:,\:\: k=2 \; ; \\[1.0 ex]
F_{k+3} + n-4 ,&      2   \, | \, n \:\:,\:\: k>2 \; .
\end{cases}
$$
\end{theorem}

\begin{corollary}
If $n \!=\! 3$ we have obtained formula
$\mbox{\textbf{\tt g}}^{\mbox{\tiny \rm G}\!}(k) = F_{k+3} + 1$ from the first section.
\end{corollary}
\begin{remark}{\rm
The values of the function $\mbox{\textbf{\tt g}}^{\mbox{\tiny \rm G}\!}(k)$ are given in {\rm \cite{Sloane}} as the following sequences
{\rm \textrm{A001611}} \mbox{$(n\!=\!3)$},  
{\rm \textrm{A000045}} \mbox{$(n\!=\!4)$},  
{\rm \textrm{A157726}} \mbox{$(n\!=\!5)$},  
{\rm \textrm{A157725}} \mbox{$(n\!=\!6)$},  
{\rm \textrm{A157729}} \mbox{$(n\!=\!7)$},  
{\rm \textrm{A157727}} \mbox{$(n\!=\!8)$},  
{\rm \textrm{A187107}} \mbox{$(n\!=\!9)$},  
{\rm \textrm{A187179}} \mbox{$(n\!=\!10)$}, 
$(k \!>\! 2)$.}
\end{remark}

\bigskip
\noindent
{\bf Acknowledgement.}
Research supported in part by the Serbian Ministry of Science, Grants No. ON174032 and III44006.

\bigskip

\break

\bigskip

\bibstyle{plain}

\medskip
\noindent
{\tt Keywords:}

\smallskip
\noindent
\mbox{\tiny $\blacksquare$} {\tt Div, Grad, Curl, Directional derivative}

\noindent
\mbox{\tiny $\blacksquare$} {\tt Differential forms, Fibonacci numbers}

\noindent
\mbox{\tiny $\blacksquare$} {\tt Compositions of the differential operations}

\bigskip
\noindent
{\tt AMS Classification:}

\smallskip
\noindent
\mbox{\tiny $\blacksquare$} {\tt Differential forms 58A10}

\noindent
\mbox{\tiny $\blacksquare$} {\tt Composition operators 47B33}

\noindent
\mbox{\tiny $\blacksquare$} {\tt Enumeration in graph theory 05C30}

\break

\bigskip
\noindent
{\tt Paper received: 22/09/2015}

\bigskip

\bigskip
\noindent
{\tt Paper accepted: 27/05/2016}

\bigskip
\noindent
{ }

\end{document}